\documentclass[12pt]{interact}

\usepackage[utf8]{inputenc}

\usepackage{array}
\usepackage{float}
\usepackage{threeparttable}
\usepackage{csquotes}
\usepackage{comment}

\usepackage{epstopdf}
\usepackage{subfigure}

\newcommand{\href}[2]{#2}

\begin{document}
\bibliographystyle{tfcad}

\title{Blended Mastery Learning in Mathematics}

\author{
\name{Timo~Pelkola\textsuperscript{a}, Antti~Rasila\textsuperscript{b} and Christopher~Sangwin\textsuperscript{c}\thanks{CONTACT C.~J.Sangwin. Email: {\tt C.J.Sangwin@ed.ac.uk}, A.~Rasila. Email: {\tt antti.rasila@iki.fi}}}
\affil{\textsuperscript{a}Department of Mathematics and Systems Analysis, Aalto University, P.O. Box 11100, FI-00076, Finland;
\textsuperscript{b}Department of Mathematics and Systems Analysis, Aalto University, P.O. Box 11100, FI-00076, Finland. +358 50 593 6016. orcid.org/0000-0003-3797-942X;
\textsuperscript{c}School of Mathematics, University of Edinburgh, The King's Buildings
Edinburgh, EH9 3FD, United Kingdom.  +44 131 650 5966. orcid.org/0000-0002-3725-8625}
}
\maketitle

\begin{abstract} 
In this paper we report a study in which we have developed a teaching cycle based closely on Bloom's Learning for Mastery (LFM).
The teaching cycle ameliorates some of the practical problems with LFM by making use of the STACK online assessment system to provide automated assessment and feedback to students.
We report a clinical trial of this teaching cycle with groups of university level engineering students.
Our results are modest, but positive: performance on the exercises predicted mastery according to the formative tests to a small extent.
Students also report being supportive of the use of the new teaching cycle.
\end{abstract}


(Word count: 5909.)

\section{Introduction}

This research is motivated by the remarkable observation of \cite{Bloom1984} that students taught by an individual tutor achieve test scores which are two standard deviations better than students who attend traditional classroom teaching.
Learning for Mastery (LFM) is an educational philosophy proposed by Bloom as a partial solution to the problem of finding resources for individual tutorials.
However, Learning for Mastery also has practical problems.
Current automatic computer aided assessment (CAA) of mathematics has reached a level of sophistication which suggests some of the practical problems with LFM might be overcome, and this is what we set out to investigate.
Can the practical problems with implementing Bloom's Learning for Mastery be overcome effectively with online CAA in mathematics?
In this paper we report a study in which we have developed a teaching cycle based closely on Bloom's LFM, making use of online assessment.
We report an action research study to investigate whether we see any significant learning gains using CAA and our LFM approach.

In Section \ref{sec:background} we provide a theoretical background to LFM and discuss contemporary CAA of mathematics in more detail.
Our precise research questions are given in Section \ref{sec:questions}.
Section \ref{sec:methodology} provides details of the methodology undertaken to address our research questions.
Results in Section \ref{sec:results} precede the final discussion.

\section{Background}\label{sec:background}

\subsection{Mathematics for university engineers}

All university engineering students learn mathematics as a core part of their undergraduate education.
Engineering mathematics curricula have been well-developed as an ongoing international collaboration, see \cite{Barry1992,Mustoe2002,Alpers2013}.
The resulting framework includes content and concepts, but goes well beyond this to include competencies.
Indeed, \cite{Alpers2013} opens the executive summary of the most recent framework document by arguing that
{\em ``the main message of this new edition is that although content remains important, knowledge should be embedded in a broader view of mathematical competencies."}
The phrase {\em ``mathematical competencies"} means that a student has proficiency in a set of interrelated mathematical skills.
The previous work of \cite[p. 116]{Kilpatrick2001}, for example, identified  conceptual understanding, procedural fluency, strategic competence, adaptive reasoning and productive disposition as five important strands.

\subsection{Mastery skills}

We separate mathematical skills (loosely) into two groups: mastery and problem solving skills (for related discussion, see \cite{Rasila2015}).
The essential distinction is that mastery skills are rarely the end goal, rather they form part of a subsequent wider task.
These skills form a loose hierarchy: weak basic conceptual and procedural skills seriously hinder a student's ability to formulate and solve mathematical problems.
\cite{Skemp1971}, for example, framed the discussion of this issue in terms of a {\em schema}: {\em ``inappropriate early schemas will make the assimilation of later ideas much more difficult, perhaps impossible}, \cite[pg 51]{Skemp1971}.
Note that mastery skills are framed within a particular context and the goals of instruction.

Mastery skills are emphatically {\em not} confined to the lower order tasks, such as recall of knowledge.
Mathematics is highly structured and one cannot jump from one level to another without doing the necessary steps between.
We include basic deductive reasoning as a mastery skill, at least to the extent that the student should understand the role of assumptions, conclusion, particular/universal statements, etc.
Without these it is impossible to create even modest chains of reasoning needed to apply more complex methods and procedures, typically taught to engineering students.
As a specific example, partial fractions require students to look ahead to anticipate the consequences of their choices.
Symbolic integration, in turn, relies on choosing particular algebraic forms, including re-writing rational terms as partial fractions.
In this context, multi-step partial fractions and symbolic integration techniques are mastery skills precisely because successful implementation of these skills are not the end point for engineers.

We should also give some idea what is not a mastery skill.
Problem solving skills are often applicable more widely, and are affective in nature (e.g.~resilience) rather than framed in terms of specific knowledge schemas.
Problem solving skills can often only be evaluated in terms of qualitative better--worse judgements, rather than right--wrong absolute judgements.
There is a substantial body of work on the learning of teaching of problem solving skills, e.g.~see \cite{Polya62}, and contemporary discussion of pedagogy for engineers \cite{Michalewicz2008}.
Since effective problem solving is normally considered to be an important part of the end goal, we do not include these skills within mastery skills.
Similarly, skills which do not form part of subsequent wider tasks are also not included within mastery skills.
Depending on the goals of the course, mastery skills may include both pen/paper calculations and the use of tools like CAS or even programming environments like MATLAB.

\subsection{Teaching, assessment and Learning for Mastery}

Different areas of mathematical proficiency require different learning strategies, e.g.~conceptual and procedural abilities are typically learned though conscious practice of exercises.
Logical thinking skills normally require significant self-reflection and discussions with a human teacher.
Similarly, advanced problem solving skills usually require a significant amount of human-to-human interaction.

Assessments, particularly high-stakes examinations, are often cited as important divers of students' learning by providing strong extrinsic motivation.
We acknowledge that high-stakes school examinations have been criticised for privileging procedural items over conceptual e.g.~\cite{Iannone2012,Noyes2011}.
At universities \cite{Tallman2016} found that little had changed in the last twenty five years:  the majority of items  required students to recall and apply a rehearsed procedure and few required conceptual understanding or problem solving.
This emphasis on procedural items is partly explained by the ease with which they can be produced and scored \cite{Swan2012}, indeed compared to other subjects scoring reliabilities tend to be high in mathematics \cite{Brooks2004}.
For further discussions of mathematical tasks see \cite{Smith1996}, \cite{Pointon02}, \cite{Watson2015} and \cite{Foster2013}.

The review of \cite{Bloom1984} considered research which compared different forms of teaching.
\cite{Bloom1984}  reports that individual tutoring resulted in student achievement which is two standard deviations better than that of students who attend traditional classroom teaching.
To close this gap \cite{Bloom1984} devised and evaluated a teaching intervention called Learning for Mastery (LFM).
In LFM students are regularly tested by using formative tests and students are required to demonstrate a correct answer to 90\% of the test problems, i.e.~demonstrate ``mastery''.
When a student falls short of mastery further teaching and testing is repeated, several times if necessary.

One of the practical impediments to LFM is the difficulty faced by the teacher who has to orchestrate the work of many students who are potentially all at different stages.
They also potentially have to devise different but related formative tests.
In traditional settings such extensive testing is still impractical.
Certainly in typical university entry-level mathematics courses, with hundreds of students, this will be impossible.
Online assessment has the potential to remove this practical barrier.
However, mastery learning can lead into surface-oriented learning strategies, especially if formative testing is mainly based on multiple choice questions.
Our interest in this topic arose because of the potential we see with contemporary online assessment in mathematics.

\begin{figure}
\centering
\resizebox*{9cm}{!}{\includegraphics{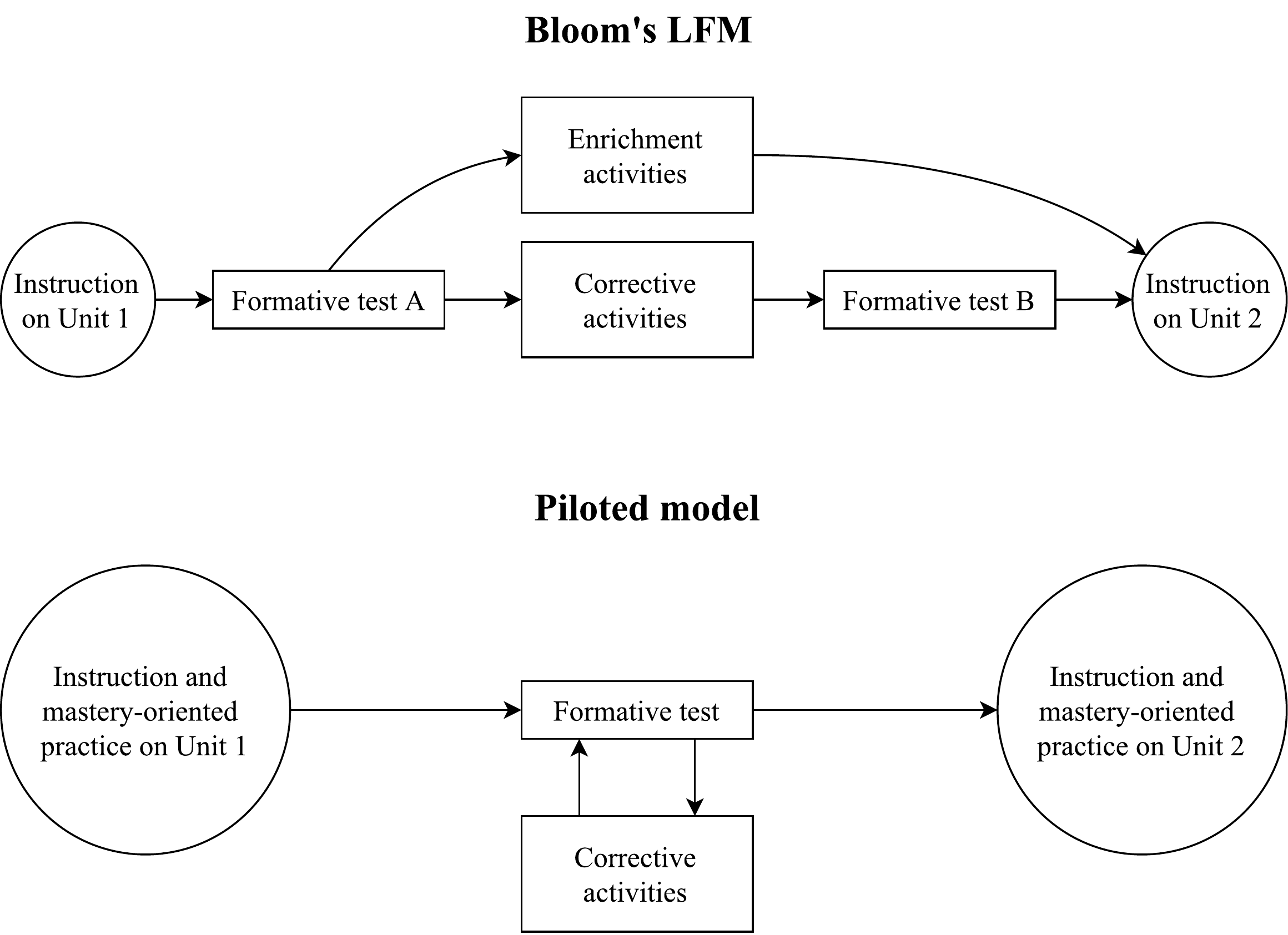}}
\caption{A comparison of Bloom's Learning for Mastery (LFM) cycle and our model} \label{fig:LFM}
\end{figure}

The current research is based on experiences gained in previous projects, such as \cite{Rasila2010}
where we started to work with the online learning system STACK as a tool for learning basic calculation techniques for engineering students, and
\cite{Majander2011} where we tried to use formative assessment (much in the sense of Bloom) to improve students' motivation to participate in the course activities.  However, these experiments lacked the corrective measures associated with Mastery Learning, which we report here. Besides improving learning outcomes, we are also interested in finding objective assessment methods suitable for distributed and distance learning (cf. \cite{Rasila2016c}).

\subsection{Online assessment for mathematics}

Computer aided assessment is well-established and widely used to support the teaching and learning of mathematics.
There is over a quarter of a century of experience developing automatic online assessment of mathematics which goes beyond relying on multiple choice and similar question types with their well-known difficulties (see \cite{2017MCQ}).
It is now common practice in online assessment of mathematics to accept a mathematical expression from a student as an answer, automatically establish relevant objective properties and provide feedback.
For example, a student might enter an algebraic expression and the teacher will have specified in advance that the computer should seek to establish algebraic equivalence with the correct answer and that it is written in a particular algebraic form, such as factored.
Normally, there are a variety of correct answers, e.g. $(x-1)(x+2)$ or $(x+2)(x-1)$.

The following features are now typical in many, if not most, mathematical systems.
\begin{itemize}
  \item Questions are randomly generated in a structured way using computer algebra systems (CAS).
        Normally the question and steps in a fully worked solution are reverse engineered from the teacher's answer.
        Quiz management components can also randomly select from a question bank to create an activity for an individual student.
  \item Students provide the final answer in the form of a mathematical expression, e.g.~an equation, rather than responding to multiple choice questions.
      It is not yet typical to automatically assess a complete argument or proof.
  \item Objective mathematical properties of answers are automatically established, e.g.~algebraic equivalence with a correct answer.
  \item Outcomes are automatically generated (including feedback) which fulfil the purposes of formative and summative assessment.
  \item Data on all attempts at one question, or by one student, are stored for later analysis.
\end{itemize}
The ability to randomly generate similar questions is particularly important for mastery learning.
Previous experience suggests the high value to students of the corresponding worked solutions, which provide a model from which students can base their answer to subsequent similar versions.

The current technical state of the art in online assessment of mathematics focuses on accepting a final answer from students and automatically establishing mathematical properties.
Many example systems provide the features we have outlines above.
An Example is the STACK online assessment system described in Section \ref{sec:stack}.  For a review of similar systems see \cite{2013CAA}.
While these systems do not (yet) fully assess complete solutions provided by students, we are aware of a number of parallel developments to implement checking of ``line by line" working in many procedural situations.
In the near future checking of line by line reasoning, and simple logic, is likely to become standard.

\section{Research questions}\label{sec:questions}

In this paper we report an action research study to investigate the following research questions.
\begin{enumerate}
\item To what extent is STACK suitable for implementing Learning for Mastery?
\item Can mastery according to formative tests be predicted from STACK exercise data?
\end{enumerate}
Lastly, we are interested in how students react to the STACK online tests used in our learning model.

\section{Methodology}\label{sec:methodology}

\subsection{Adapting Mastery learning for an online environment}

LFM suggests pairing formative assessment with appropriate correctives and we are interested in whether the practical problems with implementing LFM can be overcome effectively with online assessment.
Bloom’s Learning for Mastery model was adapted in our study using weekly online exercises and formative tests to assess mastery in core skills.
As a result, the methodology in this study differed in some ways from the original LFM model.
In LFM, mastery is assessed only with formative tests, which usually come in the form of invigilated multiple-choice questionnaires with different versions for reattempts, which are usually limited to two.
In this study mastery was assessed with online exercises.
The same formative test was used for each attempt, with the possibility of an unlimited number of attempts.
The formative test was given the name ``practice exam'' during the course, since this term was more familiar to the students.

The learning units were slightly extended to lessen the workload from the formative tests.
Also, some of the higher-order learning objectives in the course were not covered by the formative tests or online exercises, as automatic assessment of these are difficult without, in our view, fatally compromising the test validity.
Since the online component of the course covered mostly procedural skills, a new ``guided discovery'' type of project work was introduced for the exercise sessions to provide students with a balance of assessments during their course.
This consisted of four paper-based assignments and a final report about the mathematics of harmonic oscillation.

Our current study used automatic feedback generated by STACK questions which served as the primary corrective.
The formative test items were also paired with third-party videos of similar worked examples, which were made available after the first submission of the test.
We believe that students who had already gained mastery would gain some benefit from taking the formative test anyway.
As the exercises already provided some formative assessment and correctives, the necessity for separate formative tests was evaluated.

\subsection{Courses selected for the study}

Calculus I (MS-A010x) is a six-week (5 ECTS credits) compulsory course for science and engineering students at Aalto University covering single variable differential and integral calculus and ordinary differential equations. The course is offered separately for each degree programme, but with similar content.
MS-A0106 (for student majoring in mechanical and construction engineering) and MS-A0107 (for students majoring in chemical engineering)
were selected for this study, which took place as part of the continuing Aalto Online Learning (A!OLE) strategic development project coordinated at the Aalto University School of Science.
The courses consisted of four hours of lectures and exercise sessions per week, weekly online exercises, paper-based assignments, formative tests at the end of learning units and a paper-based course exam.
The course was divided into two three-week learning units, with the first unit covering limits, series and differential calculus and the second unit integral calculus and ODEs.

\subsection{Online assessment of mathematics with STACK}\label{sec:stack}

Our study adopted the STACK online assessment system.
STACK has sustained development and use for over a decade with significant contributions of code from Aalto University Finland (see \cite[Chapter 8]{2013CAA} and, for very recent work \cite{Harjula2017}), the United Kingdom Open University and latterly the University of Edinburgh in Scotland.
STACK was originally developed for Moodle but has been ported to ILIAS (see {\tt http://www.ilias.de}, retrieved May 2017) and is used in other systems, including Blackboard, through the LTI protocol.
See {\tt https://stack.maths.ed.ac.uk/demo} (retrieved May 2017).
STACK was developed by the last author, and the experimental study reported in this paper was undertaken by the first two authors at an independent institution.
The key features of STACK include its mathematical sophistication, and the full authoring interface which aims to give teachers a wide range of options in a way which still makes writing learning materials practical.

STACK is used reliably with thousands of users on over 700 registered Moodle sites.
For example, at the United Kingdom Open University during the academic year 2015-16, students attempted over 880,000 questions on seven modules.
The STACK question type accounted for approximately 15\% of all questions used, and is second only to multiple choice in popularity (at 35\% of all questions).
There are a number of large international projects such as the Abacus {\tt https://abacus.aalto.fi/} multi-lingual material bank which makes use of STACK, \cite{Rasila2016b}.
Other projects include \cite{Barbas2016}, \cite{Makela2016} and \cite{Pavia2015}, and publishers are increasingly supplementing textbooks with online assessments such as \cite{Coletta2010} which has 600 online homework problems written with STACK.

\subsection{Description of the procedure}

\begin{figure}
\centering
\resizebox*{7cm}{!}{\includegraphics{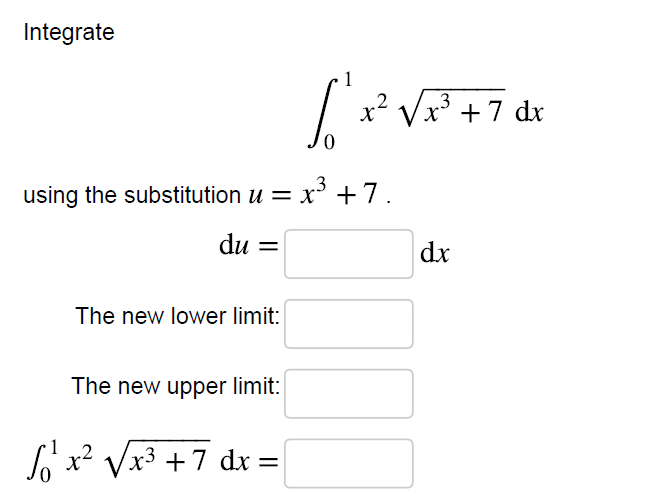}}
\caption{An example of a STACK question used on the course}
\label{fig:question}
\end{figure}

New STACK questions were developed for the formative tests and weekly online exercises.
The same set of questions were used on both courses.
The mastery threshold was set to 75-80\% of the weekly exercise or formative test points.
Neither the online exercises or formative tests were strictly compulsory, but both contributed a small part to the final course grade.
Only points above each mastery threshold were awarded, counting reattempts.
As a result, the points had minimal effect on course grading, and should be considered primarily as formative assessment.

The online exercises and formative tests had slightly different functions and were setup accordingly.
While both gave feedback on the progress of a student's learning, the online exercises were meant for initial practice, while the formative test was to ensure that mastery in those skills had actually been gained and retained.
In both cases questions could be reattempted an unlimited number of times without penalties, but in exercise questions the feedback was immediate, whereas in the formative tests it was deferred until submission of the entire test.

The weekly online exercises were due on a Sunday, and consisted of five questions related to the lectures of the week. The formative tests could be taken at the end of the learning unit before the next lecture or course exam, and the use of calculators, textbooks or other accessories were discouraged although not controlled. The first formative test included four and the second five items.

\section{Results} \label{sec:results}

95 of 134 enrolled students in MS-A0106 and 118 of 198 enrolled students in MS-A0107 consented to the use of their data for this study.
Of these students, 176 in the first learning unit and 168 in the second had opened all the weekly quizzes and the formative test at least once, which was counted as an attempt. These were used for predictive modeling.

Individual STACK item scores and numbers of attempts were extracted from the Moodle learning management system using a purpose-made export tool. This data was then imported to R for analysis.
Both the data from STACK and course feedback was used to determine the suitability of STACK for implementing mastery learning.
In the absence of a control group and proper pre/post-tests, the effectiveness of mastery learning itself was not considered in this study.

Predictive modelling was experimented with various different classification methods and pre-processing with the help of the caret R package.
Performance was evaluated with ten-fold cross validation with three repeats. Similar results were achieved with many of the methods. The results from logistic regression ('glm' in caret) are presented here.

{\em Mastery} was defined as achieving 4 out of 5 (80\%) or 3 out of 4 (75\%) points on a weekly quiz or formative test.
{\em Initial mastery} denotes the percentage of students who had achieved mastery on the first attempt, and {\em eventual mastery} those who achieved mastery on any attempt.
The difference between initial and eventual mastery is the percentage of students, who gained mastery after the first attempt.
The mastery statistics are presented in table \ref{initevent}.
On average, 88\% of quiz and 94\% of formative test takers achieved mastery eventually.
Initial mastery was achieved by 12\% on quizzes and 45\% on formative tests.
As could be expected, the level of initial mastery was higher on the formative tests than on the quizzes.
It was however significantly lower than the eventual mastery on quizzes would suggest.

\begin{table}
\tbl{Percentages of students who had gained mastery}
{\begin{tabular}{rcccccccc}
\toprule
\textbf{\underline{Initial mastery}}  & \textbf{W1}   & \textbf{W2}   & \textbf{W3}   & \textbf{FT1}  & \textbf{W4}   & \textbf{W5}   & \textbf{W6}   & \textbf{FT2}  \\
MS-A0106         & 13\% & 22\% & 7\%  & 47\% & 6\%  & 11\% & 8\%  & 52\% \\
MS-A0107         & 22\% & 25\% & 10\% & 44\% & 7\%  & 13\% & 3\%  & 38\% \\
Both            & 18\% & 24\% & 8\%  & 45\% & 6\%  & 12\% & 6\%  & 44\% \\
\textbf{\underline{Eventual mastery}} &      &      &      &      &      &      &      &      \\
MS-A0106         & 92\% & 95\% & 85\% & 95\% & 79\% & 90\% & 90\% & 98\% \\
MS-A0107         & 90\% & 90\% & 83\% & 90\% & 85\% & 90\% & 87\% & 94\% \\
Both            & 91\% & 92\% & 84\% & 92\% & 82\% & 90\% & 89\% & 96\% \\
\textbf{\underline{Difference}}       &      &      &      &      &      &      &      &      \\
MS-A0106         & 80\% & 73\% & 78\% & 47\% & 73\% & 78\% & 82\% & 46\% \\
MS-A0107         & 67\% & 65\% & 74\% & 47\% & 78\% & 77\% & 84\% & 56\% \\
Both            & 73\% & 69\% & 76\% & 47\% & 76\% & 77\% & 83\% & 52\% \\
\bottomrule
\end{tabular}}
\label{initevent}
\end{table}

\begin{figure}
\centering
\resizebox*{13cm}{!}{\includegraphics{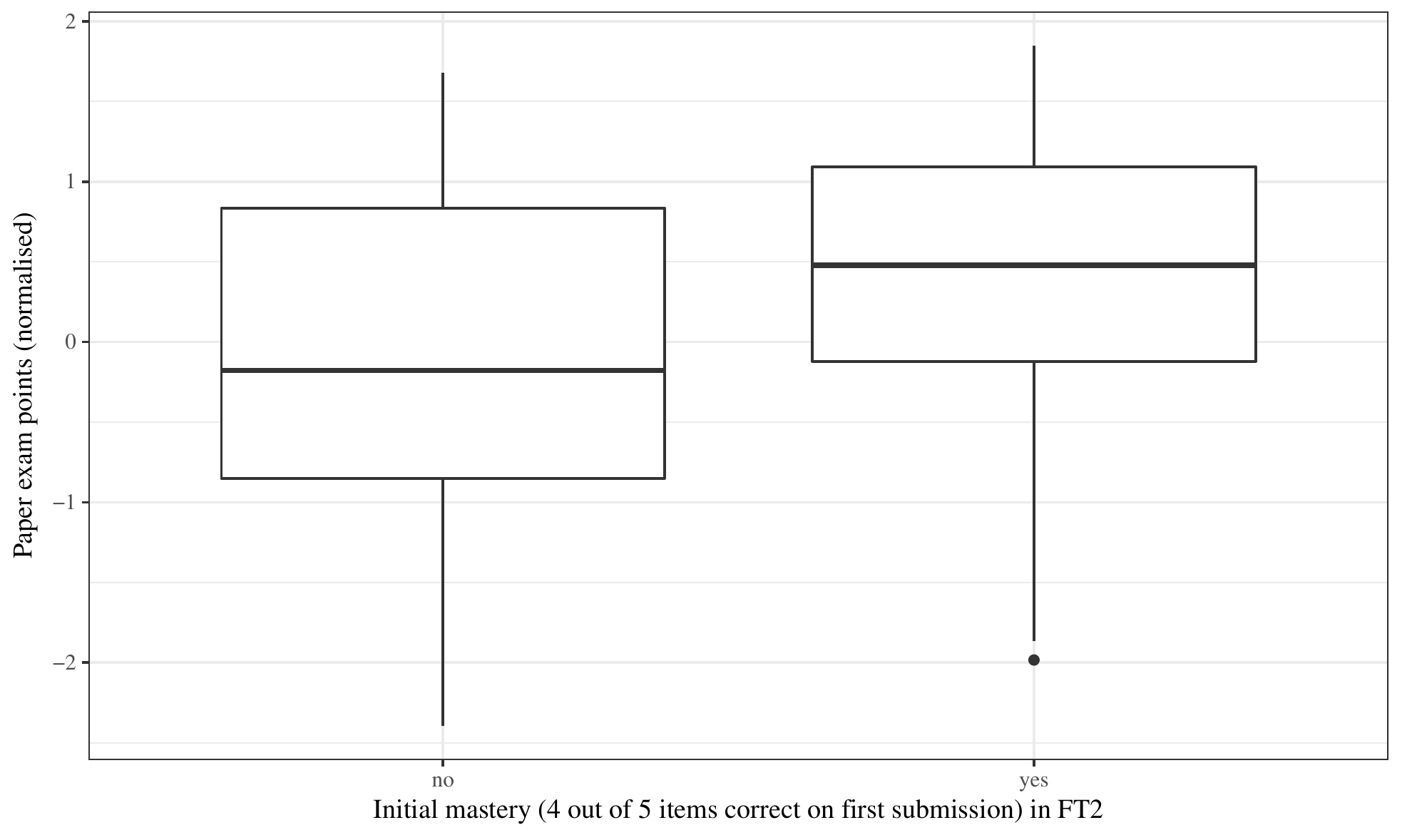}}
\caption{Examination points compared between mastery and non-mastery students}
\label{fig:examplot}
\end{figure}

When the pen-and-paper examination scores were compared against initial mastery on the second formative test (Figure \ref{fig:examplot}), a difference of 0.51 standard deviation in mean test scores was found. This would further suggest that eventual mastery is not entirely equivalent to initial mastery.

\subsection{Qualitative questionnaires}

The MS-A0106 course feedback questionnaire included four likert-scale questions considering the ML model used on the course, and a summary of results is shown in Table \ref{table:questionnaire}.
The feedback was mostly in favour of the model.

The mastery bonus point scheme (no points awarded below mastery) seemed to encourage (40\% of respondents) more than discourage (14\%) practice. 37\% found the formative tests very useful, while only 5\% found the formative tests not useful at all.

The videos that served as correctives on the formative tests were also found useful (91\%) by those who had watched the videos (26\%). It is unclear why so many chose not to watch the videos, but the figure should be nonetheless compared to the level of initial mastery (50\% on average in MS-A0106), as they were only watchable after the first attempt.

As the formative tests were meant to be solved without the aid of calculators or learning materials, but were not invigilated, activities which might be considered as `cheating' caused some concern.
A majority (84\%) of students admitted using accessories like calculators and books during the formative test at least a few times, although only 13\% reported this often.
Judging from the mass of erroneous answers even to the questions easily solvable with a CAS, it would seem that at least the first attempt was usually relatively sincere.

\begin{table}
\tbl{Mastery learning -related questions from the course feedback questionnaire (89 respondents)}
{\begin{tabular}{m{20em} c c c c c}
\toprule
& \textbf{1}    & \textbf{2}    & \textbf{3}    & \textbf{M}   & \textbf{SD}  \\ \hline
\textbf{Were the practice exams useful?}\newline 1) not at all 2) somewhat 3) very                                                                                 & 5\%  & 58\% & 37\% & 2.3 & 0.6 \\
\textbf{Were the practice exam related videos useful?\textsuperscript{*}
} \newline1) not at all 2) somewhat 3) very                                                                  & 9\%  & 41\% & 50\% & 2.4 & 0.7 \\
\textbf{Did you use accessories (calculators, books etc...) in the practice exam?} \newline1) never 2) a few times 3) often                                        & 15\% & 71\% & 13\% & 2.0   & 0.5 \\
\textbf{Mastery bonus point scheme (0 points if less than 80\% done) had mostly … to my practice} \newline 1) a negative effect 2) no effect 3) a positive effect & 14\% & 47\% & 40\% & 2.3 & 0.7 \\
\bottomrule
\end{tabular}}
    \begin{tabnote}
    \textsuperscript{*} including only those who reported watching the videos (26\% of formative test takers)
    \end{tabnote}
\label{table:questionnaire}
\end{table}

In the responses to the question \emph{``Which things were good on the course? What promoted your learning?''} parts of the LFM model were commended.
Almost all of the 74 responses mentioned exercises or exercise sessions in some way.
18 mentioned STACK exercises specifically and 8 the formative tests.
Some examples (translated from Finnish to English by us) were:

\begin{displayquote}
\emph{STACK exercises and practice exams were a good addition. Altogether all kinds of extra homework helps, since in my case drilling the basics should be emphaised a bit more before moving on to applications.}
\newline\newline
\emph{The middle exams gave a good sense of how well you have mastered the course content.}
\newline\newline
\emph{The practice exams forced [me] to revise.}
\end{displayquote}

Also the mastery-oriented bonus point scheme got mentioned:

\begin{displayquote}
\emph{A good thing on the course was that the STACK exercises were, in a way, mandatory.}
\end{displayquote}

There was also a counterpart to the previous question (\emph{Which things were bad / didn’t work? What hindered your learning?}). The 68 responses were mostly focused on the project assignments, lectures and lecture notes. Two students felt there were too many different types of activities on the course.

STACK exercises were mentioned to be both too difficult and not challenging enough.

\begin{displayquote}
\emph{... Also some of the STACK exercises were such that I couldn't find even a hint of a “basic exercise” in those. At least the lectures gave me no clue of solution models, and sometimes I didn't get it even after the teaching assistant had explained it.}
\newline\newline
\emph{... There were all too many exercises and they all were unchallenging. I’d prefer three times less exercises but more challenging ones. Especially STACK exercises often felt like a waste of time.}
\end{displayquote}

The formative tests were not criticised apart from unclear instructions.

\subsection{Predictive Modelling}
Predicting mastery on the formative tests based on prior performance on the quizzes proved to be more challenging than anticipated.

\begin{figure}
\centering
\resizebox*{13cm}{!}{\includegraphics{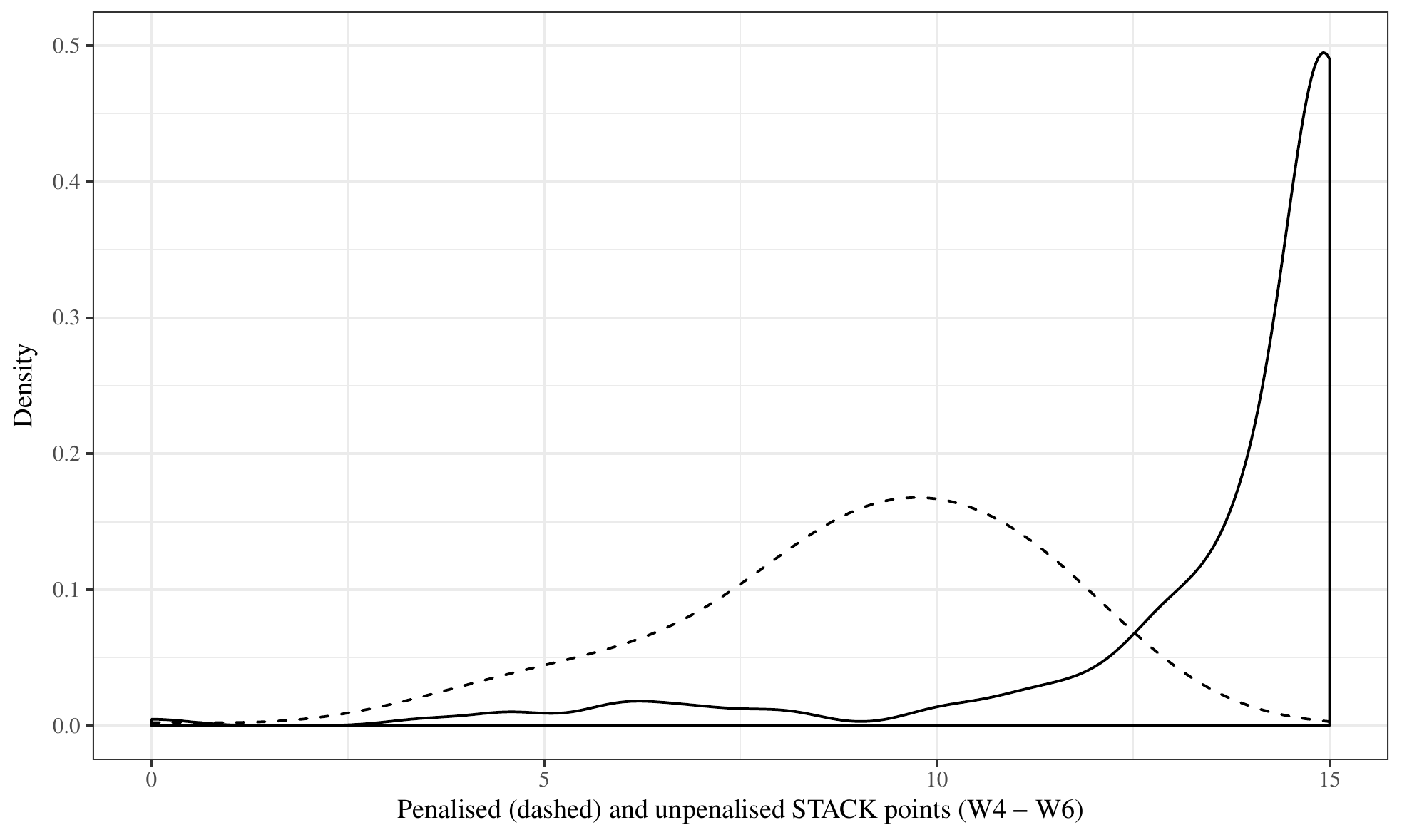}}
\caption{Ceiling effect on unpenalised points from weeks 4 to 6}
\label{fig:distribution}
\end{figure}

A notable ceiling effect was observed with the unpenalised quiz points. Simulated penalty was later applied with a formula
\begin{equation}\label{penal}
\text{penalised~points} = \text{floor}\left(\text{raw~points}\right) \cdot 0.7^{\text{reattempts}},
\end{equation}
where $\text{floor}(x)$ rounds partial points down towards zero. The formula resulted in a less skewed distribution, shown in Figure \ref{fig:distribution}. The penalised points from quizzes 4-6 also had a higher Spearman correlation (0.51) with the paper exam points than the unpenalised raw points (0.40).

It should be noted that actual penalties, or limiting the number of attempts, are likely to have some effect on behaviour.
High numbers of attempts were observed in some cases, suggesting that some students adopt a trial-and-error strategy when such behaviour goes unpenalised.

When comparing students’ exercise point sums against initial mastery on the formative tests, it could be seen that the points provided poor separation between mastery and non-mastery. The difference in median points were highest when one reattempt was allowed, but the ceiling effect became apparent with further attempts.

\begin{figure}
\centering
\resizebox*{13cm}{!}{\includegraphics{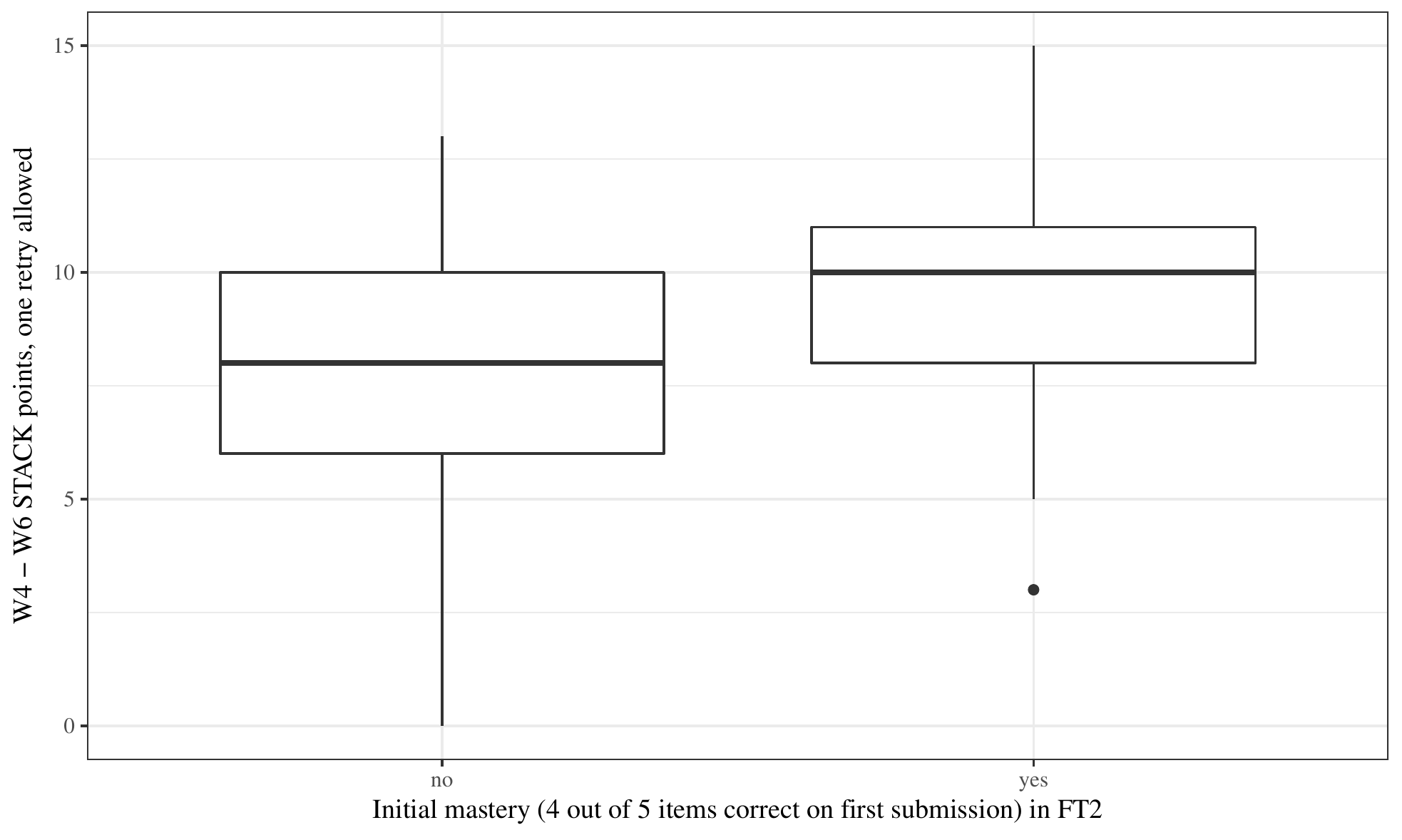}}
\caption{When one reattempt per exercise problem was taken into account, exercise points between mastery and non-mastery students provided some  separation}
\label{fig:bp1}
\end{figure}
\begin{figure}
\centering
\resizebox*{13cm}{!}{\includegraphics{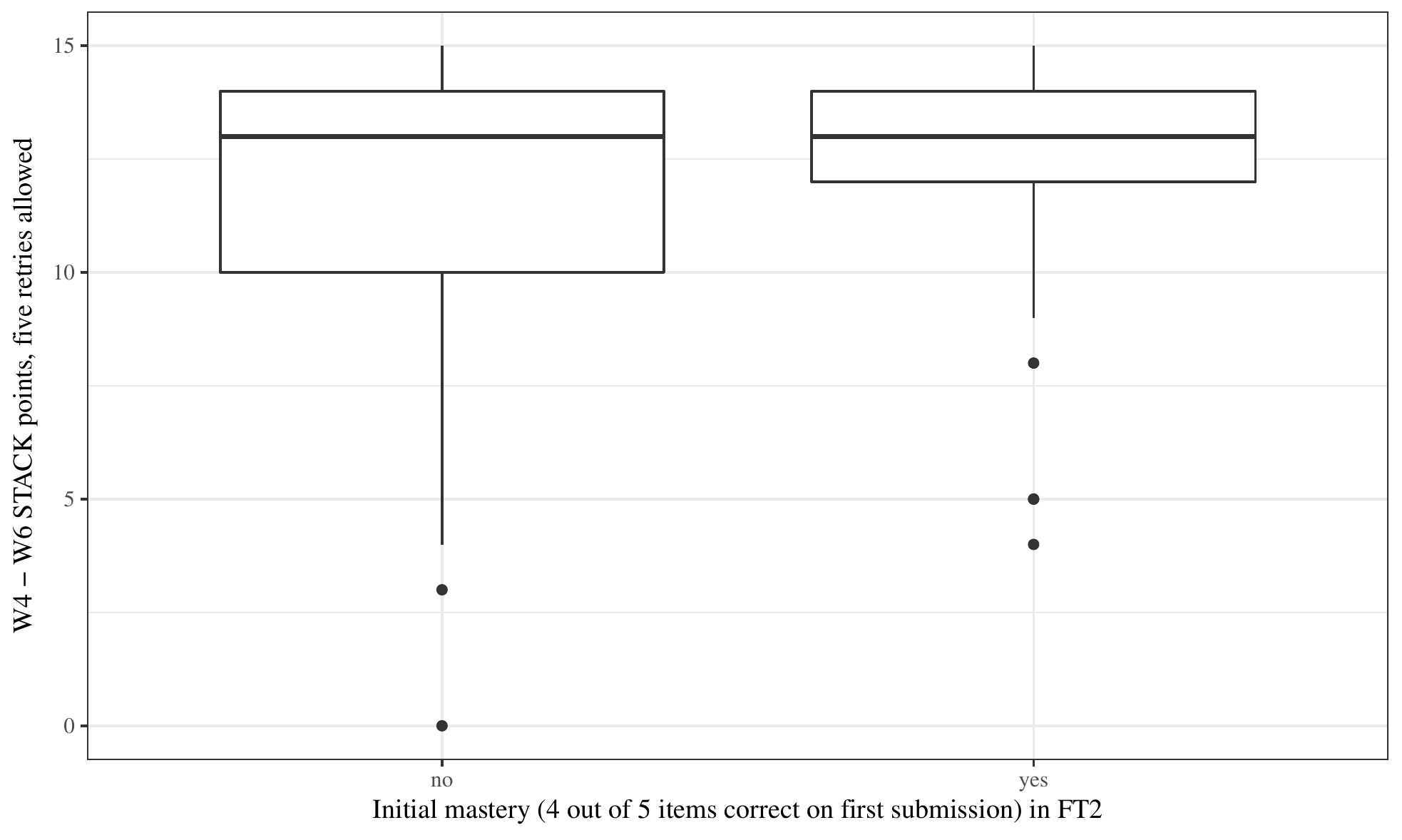}}
\caption{After five reattempts there is no more difference in median points}
\label{fig:bp5}
\end{figure}

As could be concluded from the data in table \ref{initevent}, the eventual mastery in the quizzes did not translate into initial mastery in the formative tests, and the sum of points did not seem to separate mastery and non-mastery either (Figures \ref{fig:bp1} and \ref{fig:bp5}).
Therefore, a more sophisticated model would be required to tell whether a student would be likely to achieve mastery in the following formative test.
An attempt was made to construct a unit mastery classifier that could ultimately replace the formative tests.

We experimented using various different methods found in the caret R package.
Logistic regression performed comparably to some of the more advanced methods such as gradient-boosted trees and was chosen for the model.

Logistic regression has the additional advantage of providing class probabilities, which allows us to optimise the classification threshold easily. In this case, the cost of a false positive (inadequate learning) could be considered greater than that of a false negative (waste of time).

Data from quizzes 4-6 were used to predict the initial mastery on the second formative test. In the end, the sum of penalised points (equation \ref{penal}) provided the best results.

It should be noted that the number of complete observations (168 in the second learning unit) limits how many predictors can be used without overfitting, and might have been the reason why the individual question points and numbers of reattempts did not result in a more accurate model.

Some pre-processing of the data was also needed, because the number of reattempts before success and giving up are measuring essentially different things. The exponential penalty scheme (equation \ref{penal}) was chosen after some experimentation, as this provided a way of reducing points and number of reattempts into a single variable and did not suffer from a floor effect as would a linear model.

The resulting model, predicting that a student would not achieve initial mastery on the second formative test, had an accuracy of 0.64 which is a small improvement over predicting that no student would achieve mastery (0.56).

\begin{table}
    \tbl{Confusion matrix of the classifier (10-fold cross validation with 3 repeats)}
    {\begin{tabular}{r c c}
    \toprule
        &\multicolumn{2}{c}{\textbf{\underline{Actual}}}\\
        \textbf{\underline{Prediction}}&non-mastery&mastery\\

        non-mastery&33.7\%&17.1\%\\
        mastery&18.7\%&30.6\%\\
        \bottomrule
    \end{tabular}}
    \label{table:confusion}
\end{table}

\section{Conclusion}

\subsection*{Is STACK suitable for implementing mastery learning?}

The STACK system is able to assess most of the learning objectives of Calculus I, and as such is in theory suited for implementing ML on the course. From a technical perspective STACK has many advantages over other similar online assessment systems, particularly in the potential to create sophisticated feedback.
However, we believe any online assessment system accepting algebraic answers as students' answers is likely to generate similar overall results.

The implementation was also proven to work in practice, since on each formative test and weekly quiz a considerable portion of students' achievement was raised from non-mastery to mastery (69-83\% on the quizzes and 47-52\% on the formative tests).

Based on the course feedback, students generally approved of the model.
The formative tests were seen as useful and the mastery-oriented bonus point scheme encouraged the students as was intended.
However, some concern is caused by the fact that eventual mastery on the weekly quizzes did not translate to initial mastery on the formative tests, and that those who had achieved initial mastery on the second formative test also did better on the paper examination. This could be due to a difference between \emph{exercise} and \emph{test proficiency}.

Solving an exercise problem might be considerably easier than solving the same problem in a test situation for a number of reasons. The student may get help from a peer or a teacher, does not have to rely only on memorised facts, can check his answer, reattempt and may also be more inclined to use a calculator.
Similarly, reattempts of a test may also be fundamentally different from the first attempt.

Even so, there is no definite answer to which one of these is the desired level of proficiency.
The formative tests however do seem to reveal something the exercises alone cannot, and thus could be beneficial to learning in any case.

The difference between initial and eventual mastery could also be blamed on the ineffectiveness of the correctives or the fact that eventual mastery may have been achieved with the aid of a calculator.
Paper-based examinations have been refined for many centuries, but using online assessment effectively is in its infancy.

\subsection*{Can mastery according to formative tests be predicted from STACK exercise data?}

Performance on the exercises predicted mastery according to the formative tests to a small extent, and in this case does not warrant using a predictive model as a replacement for the formative tests.
However, the result was still positive and could possibly be further improved with more observations, different independent variables and fine-tuning of the model.
Some of the considerations from the previous section also applies here. Invigilation of the formative tests could make the model training data more reliable.

Our results also suggest that STACK-based examinations are not a completely realistic substitute for pen and paper examinations, which is a rather significant result.

\section*{Acknowledgements}

The authors acknowledge the financial support of their departments in for this work.

\section*{Disclosure statement}

There is no conflict of interest reported by the authors.

\section*{Biographical note}

{\bf T.~Pelkola}

Timo Pelkola is a Research Assistant working at Aalto University, Helsinki, Finland, in the research group of A. Rasila. He is presently a Master Student at University of Helsinki, from where he is expected to graduate as a Mathematics and Computer Science teacher in the year 2017. His topics of interest include e-assessment in mathematics, and use of interactive technologies in supporting student motivation.

{\bf A.~Rasila}

Antti Rasila is a mathematician working at Aalto University, Helsinki, Finland as a Senior University Lecturer and the leader of the Research Group of Computer Aided Mathematics Teaching.  He is also a Docent of mathematics at University of Vaasa and a Guest Professor at Hengyang Normal University. His fields of interest include computer-aided methods in mathematics education, but also certain topic in pure and applied mathematics like complex analysis, partial differential equations and numerical conformal mappings.

{\bf C.~J.~Sangwin}

Chris Sangwin is Professor in the School of Mathematics at the University of Edinburgh, United Kingdom.
His learning and teaching interests include (i) automatic assessment of mathematics using computer algebra, and (ii) problem solving using Moore method and similar student-centred approaches.

\end{document}